\documentclass[leqno,12pt]{article}
\begin{document}

\newcommand {\s}[1]{{\cal #1}} \newcounter{sectiune}
\setcounter{sectiune}{-1}
\newcounter{lemma}
\setcounter{lemma}{0}
\newenvironment{sectiune}{\setcounter{lemma}{0}%
                         \setcounter{equation}{0}%
                         \stepcounter{sectiune}%
                         \bigskip \noindent
                         \bf%
                         \arabic{sectiune}.  }%
                         \bigskip

\renewcommand{\theequation}{\arabic{sectiune}.\arabic{equation}}

\newenvironment{lemma}{
                      \stepcounter{lemma}\begingroup \noindent
                        \bf
                       (\arabic{sectiune}.\arabic{lemma})
                       Lemma.\it}{
                        \endgroup\vskip\baselineskip}

\newenvironment{proposition}{
                      \stepcounter{lemma}\begingroup \noindent
                        \bf
                       (\arabic{sectiune}.\arabic{lemma})
                       Proposition. \it}{
                        \endgroup\vskip\baselineskip}

\newenvironment{corollary}{
                       \stepcounter{lemma}\begingroup \noindent
                         \bf
                       (\arabic{sectiune}.\arabic{lemma})
                       Corollary. \it}{
                        \endgroup\vskip\baselineskip}

\newenvironment{theorem}{
                      \stepcounter{lemma}\begingroup \noindent
                        \bf
                       (\arabic{sectiune}.\arabic{lemma})
                       Theorem. \it}{
                        \endgroup\vskip\baselineskip}

\newenvironment{remark}{
                      \stepcounter{lemma}\begingroup \noindent
                      \bf
                       (\arabic{sectiune}.\arabic{lemma})
                       Remark. \rm}{
                        \endgroup\vskip\baselineskip}



\noindent
{\large \bf On the vanishing of higher syzygies of curves}\\\\
{\bf Marian Aprodu}\\\\
{\small Romanian Academy, Institute of Mathematics "Simion Stoilow",
P.O.Box 1-764, RO-70700, Bucharest, Romania (e-mail:
Marian.Aprodu\char64 imar.ro) \&\\
Universit\'e de Grenoble 1,
Laboratoire de Math\'ematiques, 
Institut Fourier BP 74,
38402 Saint Martin d'H\`eres Cedex,
France (e-mail: aprodu\char64 mozart.ujf-grenoble.fr)}



\begin{sectiune}
Introduction, main results
\end{sectiune}

\noindent
A main challenge in the theory of syzygies is to
interpret the information carried by the graded
Betti numbers of a smooth projective variety.
Notably, the attempt to understand the way
in which the distribution of zeroes in a
Betti table interacts with the geometry of the variety
has led to a considerable amount of work, motivated by the conjectures
that Green, Green-Lazarsfeld, and others had formulated
(see, for example, \cite{Gr1}, \cite{GL1}, \cite{EL},
\cite{La2}).

One of the most significant conjectures made by Green and Lazarsfeld
(cf. \cite{GL1} 3.7; see also \cite{Gr3} 3.5 and
\cite{La2} 2.3), nowadays known as the {\em gonality conjecture},
predicts that one could read off the gonality of a
smooth complex projective curve from the minimal
resolution of any line bundle of sufficiently large degree.
In order to give a precise statement, the authors
introduced the vanishing property  $(M_k)$ (see
\cite{GL1}, \cite{Gr3}), which is the following
(we use the notation of Section 1).
\\\\
{\bf Definition.} (Green-Lazarsfeld)
If $X$ is a smooth complex projective curve of genus $g$,
$L$ a line bundle on $X$, and $k\geq 0$ an integer,
one says that $(X,L)$ {\em satisfies the property}
$(M_k)$ (or, simply, $L$ {\em satisfies} $(M_k)$)
if $K_{p,1}(X,L)=0$, for all $p\geq h^0L-k-1$.

\bigskip

It is well-known (see, for example, \cite{GL1}, \cite{GL2}
or \cite{Sch}) that if $X$ carries a $g^1_k$, then
no line bundle of sufficiently large degree
can satisfy $(M_k)$. The gonality conjecture states a converse
of this fact: if $\mbox{deg}(L)>>2g$, and $(M_k)$
fails for $L$, then $X$ carries a $g^1_q$, with $q\leq k$.
Green has shown this conjecture holds for
$k=1,2$ (cf. \cite{Gr1}); Ehbauer proved it in the case $k=3$
(cf. \cite{Ehb}). Thus hyperelliptic and trigonal curves are
characterized by syzygies.

The purpose of the present work is to prove that,
under certain conditions, the property
$(M_k)$ is preserved when we add an effective divisor
to a line bundle. The first result is:

\bigskip


\noindent
{\bf Theorem 1.}
{\em Let $X$ be a smooth complex projective
curve of genus $g\geq 1$, $L_0$ be a nonspecial
globally generated line bundle
on $X$, and $k\geq 0$ be an integer such that the pair
$(X,L_0)$ satisfies the property
$(M_k)$. Then, for any effective divisor
$D$ on $X$, the pair $(X,L_0+D)$ satisfies the property $(M_k)$.
In particular, for any line bundle
$L$ with $\mbox{\rm deg}(L)\geq \mbox{\rm deg}(L_0)+g$,
the pair $(X,L)$ satisfies the property $(M_k)$ as well.}

\bigskip

There are a number of immediate consequences
of Theorem 1. For instance, if $X$ carries
a $g^1_k$, and $L_0\in\mbox{Pic}(X)$ is nonspecial,
and globally generated, then $(X,L_0)$ cannot
satisfy $(M_k)$. Further, the least $k$
for which the property $(M_k)$ fails for
a nonspecial, globally generated line bundle
cannot decrease when adding effective divisors,
and cannot pass over the gonality of $X$,
and thus it must be constant when the degree of the
line bundle grows
large enough. The challenge of the gonality conjecture
is to show that this constant always equals 
the gonality of $X$.

In view of Theorem 1,
verifying the gonality conjecture for a given $k$-gonal
$X$ of genus $g$ reduces to finding a single nonspecial,
globally generated line bundle
(for example, a line bundle whose degree is sufficiently large
compared to $2g$), with property  $(M_{k-1})$. Concretely, we get
the following criterion for testing the gonality conjecture:
\\\\
{\bf Corollary 2.} {\em  Let $X$ be
a smooth complex projective curve of genus $g\geq 1$, which
carries a $g^1_k$, and $L_0$ be a nonspecial,
globally generated line bundle on $X$ satisfying the property $(M_{k-1})$.
Then $X$ is $k$-gonal, and the gonality conjecture
is valid for $X$.}

\bigskip

By semicontinuity of minimal resolutions, and irreducibility of the moduli
space ${\cal M}_{g,k}$ of $k$-gonal curves of genus $g$
(cf. \cite{F}), it turns out that the gonality conjecture
would also be valid for a generic curve in ${\cal M}_{g,k}$,
once it was verified for a particular $k$-gonal curve of genus $g$.
It seems thus very reasonable to predict that for any positive integers
$k$, and $g$, such that ${\cal M}_{g,k}\not=\emptyset$,
a generic $k$-gonal curve of genus $g$ verifies
the gonality conjecture (see below).

One could address now the question of what happens if
we drop the nonspeciality condition in Theorem 1.
A partial answer is given by:
\\\\
{\bf Theorem 3.}
{\em Let $X$ be a smooth complex projective
curve of genus $g\geq 1$, and $k\geq 0$ be an integer
such that the pair $(X,K_X)$ satisfies the property
$(M_k)$. Then, for any effective divisor
$D$ on $X$, the pair $(X,K_X+D)$ satisfies the property $(M_k)$.}

\bigskip

Unfortunately, the result just stated is weaker than it might look
like at a first sight - it does not show that Green's generic
canonical conjecture
implies the gonality conjecture for a generic curve, as one could think of. 
The following better version of it would do the job instead.
\\\\
{\bf Conjecture.} {\em Let $X$ be a smooth complex projective
curve of Clifford dimension one, and $p$ an integer such that
$K_{p,1}(X,K_X)=0$. Then there exist a positive integer $d$, and
an effective divisor $D$ on $X$, of degree $d$, such that 
$K_{p+d-1,1}(X,K_X+D)=0$.}

\bigskip

This conjecture is obviously false if we drop the condition
that the Clifford dimension be equal to one, as seen by analyzing the case
of smooth plane curves. In exchange, it holds for other curves which are
not plane curves, such as curves lying on Hirzebruch surfaces
(see Sections 6 and 8), and, more generally, it is true for curves
which verify both the gonality conjecture, and
Green's canonical conjecture (trigonal curves, for instance). 

\bigskip

The main ingredient we use to prove Theorem 1, and Theorem 3,
is projection of syzygies, concept which was introduced
by Ehbauer \cite{Ehb}, in a coordinate-based manner.
In Section 2, we propose a 
more abstract view on the subject, in the spirit
of \cite{Gr2} 1.b.1, and we think of  projections of syzygies
as being {\em corestrictions} of the fiber-restrictions of a
certain morphism between vector bundles. We investigate
some properties of the projection morphisms, and
we show, among other things, that, under
some conditions which are almost always satisfied,
any nonzero syzygy survives by projection from
a generic point (see (2.5)).

The third Section deals with projection of
syzygies of varieties. We recall here Ehbauer's
approach, and examine the case of syzygies of curves.

In Section 4 we complete the proofs
of Theorem 1, and Theorem 3, and further,
in the final part of the paper, we analyze
some very concrete cases.

We verify first the gonality conjecture for smooth plane curves, and
for smooth curves lying on a Hirzebruch surface, in which cases
we can make use of the geometry of the ambient surface
to produce suitable line bundles.
Additionally, we recover the description of
the minimal pencils in these cases, which has
already been known before (for plane curve, see, for
example \cite{ACGH}, for curves on a Hirzebruch
surface we refer to \cite{Ma}).

In the fifth Section, we test the gonality conjecture for
nodal curves on ${\bf P}^1\times {\bf P}^1$, case which
eventually shows the following
(compare with the main result of \cite{Sch2}):
\\\\
{\bf Theorem 4.} {\em
For any integer $k\geq 3$, the gonality conjecture
is valid for a generic $k$-gonal curve of genus
$g>(k-1)(k-2)$. 
}

\bigskip

We conclude this paper by applying the vanishing result (6.2)
proved here to show that Green's canonical conjecture 
holds for smooth curves on Hirzebruch surfaces (compare with \cite{Lo}).
As a general philosophy, we expect that smart choices of line bundles on
a surface, whose restrictions  satisfy the vanishing property required in
the gonality conjecture for some curves lying on that surface, be used to
prove Green's canonical conjecture for such curves.


\begin{sectiune}
Some notation,  preliminaries
\end{sectiune}

\noindent For many of the theoretical facts included in this
section we refer to the papers \cite{Gr1}, \cite{Gr3}, \cite{La1},
without further mention. We start with $V$ a
finite-dimensional complex vector space, we denote $SV$ the
symmetric algebra of $V$, and we consider
$B=\mathop\bigoplus\limits_{q\in {\bf Z}}B_q$ a graded
$SV$-module. Then there is a naturally defined
complex of vector spaces, called the {\em Koszul complex} of $B$,
$$
...\longrightarrow B_{q-1}\otimes
\bigwedge ^{p+1}V\stackrel{d_{p+1,q-1}}{\longrightarrow}
B_q\otimes \bigwedge ^pV
\stackrel{d_{p,q}}{\longrightarrow}
B_{q+1}\otimes\bigwedge ^{p-1}V
\longrightarrow ... ,$$
whose cohomology is denoted by 
$$ K_{p,q}(B,V)=\mbox{Ker }d_{p,q}/ 
\mbox{Im }d_{p+1,q-1}. $$ 
The dimensions of the Koszul cohomology spaces
are called {\em graded Betti numbers of $B$}, 
and their elements are called {\em syzygies}.

An important property of Koszul cohomology is its functorial behavior -
any morphism of graded $SV$-modules, $A\rightarrow B$, 
canonically induces linear maps $K_{p,q}(A,V)\rightarrow 
K_{p,q}(B,V)$. Moreover,
to any short exact sequence of graded $SV$-modules,
$0\rightarrow A\rightarrow B
\rightarrow C\rightarrow 0$
one associates a long exact sequence, for any $p$
(cf. \cite{Gr1}, 1.d.4):
$$
...\rightarrow K_{p+1,0}(C,V)
\rightarrow K_{p,1}(A,V)
\rightarrow K_{p,1}(B,V)\rightarrow K_{p,1}(C,V)
\rightarrow K_{p-1,2}(A,V)\rightarrow ...
$$
\begin{remark} With the notation above,
we see that if $C_0=0$, then the map $K_{p,1}(A,V)
\rightarrow K_{p,1}(B,V)$ is injective, and if $C_0$, and
$C_1$ both vanish, then $K_{p,1}(A,V)\cong K_{p,1}(B,V)$.
\end{remark}

\noindent
{\bf Notation.} If $X$ is an irreducible complex projective variety,
$L\in \mbox{Pic}(X)$ is a line
bundle, $\cal F$ is a coherent sheaf,
and $V\subset H^0(X,L)$, 
take  $B=\bigoplus\limits_{q\in{\bf Z}} H^0(X,{\cal F}\otimes qL)$
(here, $qL$ is the $q$-th tensor power of $L$ in $\mbox{Pic}(X)$),
and denote $K_{p,q}(X,{\cal F},L,V) =K_{p,q}(B,V)$. 
If $V=H^0(X,L)$, we drop $V$ and write
$K_{p,q}(X,{\cal F},L)$, if ${\cal F}\cong {\cal O}_X$, we
suppress it, and write $K_{p,q}(X,L,V)$; the notation
$K_{p,q}(X,L)$ corresponds to the choice $V=H^0(X,L)$ and ${\cal
F}\cong {\cal O}_X$.

\bigskip

\begin{remark} 
If $L$ is a globally generated line bundle over the
smooth irreducible variety $X$, and $V=H^0(X,L)$, by denoting $X'$
the image of $X$ in ${\bf P}V^*$, through the morphism given by
the {\em complete} linear system $|L|$, and $S_{X'}=\mbox{Im}\big(
SV\rightarrow \bigoplus H^0(X,qL)\big)$ its homogeneous coordinate
ring, then, in virtue of $(1.1)$, we have natural isomorphisms
$K_{p,1}(S_{X'},V) \cong K_{p,1}(X,L)$, for all integers $p$.
\end{remark}
\begin{remark}
If $X$ is a smooth irreducible projective variety
with $h^1{\cal O}_X=0$, $L\in \mbox{Pic}(X)$
a line bundle, and $Y\in |L|$ is irreducible,
then $K_{p,1}(X,L)\cong K_{p,1}(Y,L|_Y)$
for any integer $p$. For the proof, we apply 
(1.1), and argue as in 3.b.7 of \cite{Gr1}.
\end{remark}

\noindent
{\bf Convention.} Unless otherwise stated, a {\em curve}
will always mean a smooth, connected, complex, projective
curve.


\begin{sectiune}
Projections of syzygies at large
\end{sectiune}

\noindent
We consider, as in the previous section,
$V$ a finite-dimensional complex vector space, and $B$ a graded
$SV$-module. We denote by ${\bf P}={\bf P}V^*=\mbox{Proj}(SV)$, and 
$L={\cal O}_{\bf P}(1)$.
A point $x\in {\bf P}$, which corresponds to a short
exact sequence of vector spaces:
$$
0\longrightarrow W_x
\longrightarrow V \stackrel{\mbox{ev}_x}{\longrightarrow}
L_x\longrightarrow 0,
$$
naturally induces, for any integer $l$,
a short exact sequence of Koszul complexes:
$$
0
\longrightarrow
B_*\otimes \bigwedge^{l+1-*}W_x
\longrightarrow
B_*\otimes\bigwedge^{l+1-*}V
\longrightarrow
L_x\otimes B_*\otimes \bigwedge^{l-*}W_x
\longrightarrow
0,
$$
which gives rise to a long exact sequence
(compare with \cite{Gr2} 1.b.1):
$$
(*)...\rightarrow
K_{p+1,q}(B,W_x)\stackrel
{\eta_{x}}{\rightarrow}
K_{p+1,q}(B,V)\stackrel
{\pi_{x}}{\rightarrow}
L_x\otimes K_{p,q}(B,W_x)\stackrel
{\mu_{x}}{\rightarrow}
K_{p,q+1}(B,W_x)\rightarrow
...,
$$
where $p+q=l$. We call the map $K_{p+1,q}(B,V)\stackrel{\pi_{x}}
{\longrightarrow} L_x\otimes K_{p,q}(B,W_x)$
{\em projection of syzygies centered in $x$},
and the elements of
its image, {\em projected syzygies}.

Projection of syzygies is functorial, that
is, for any morphism $B\rightarrow C$ of graded $SV$-modules,
the induced morphisms between syzygies
$K_{p+1,q}(B,V)\rightarrow K_{p+1,q}(C,V)$, and
$L_x\otimes K_{p,q}(B,W_x)\rightarrow
L_x\otimes K_{p,q}(C,W_x)$
commute with the corresponding projection morphisms.
Furthermore, projection of syzygies is compatible
with the connecting morphisms arising from the
long cohomology sequences: if
$0\rightarrow A\rightarrow B\rightarrow C\rightarrow 0$
is an exact sequence of  graded $SV$-modules, then the connecting
morphisms $K_{p+1,q}(C,V)\rightarrow K_{p,q+1}(A,V)$,
and $L_x\otimes K_{p,q}(C,W_x)\rightarrow
L_x\otimes K_{p-1,q+1}(A,W_x)$ commute with projections.

\bigskip

\begin{remark}
The composed map $\partial_x:=\big( \eta_x\otimes
\mbox{id}_{L_x}\big) \circ\pi_x:
K_{p+1,q}(B,V) \rightarrow L_x\otimes K_{p,q}(B,V)$ is the
fiber-restriction over $x$ of a natural sheaf morphism
between two vector bundles:
$$
{\cal O}_{\bf P}\otimes
K_{p+1,q}(B,V)
\stackrel{\partial}{\longrightarrow}
{\cal O}_{\bf P}(1)\otimes K_{p,q}(B,V).
$$
In order to prove this, we use the Koszul complex of $B$, 
and suitable shifts of the well-known exact complex 
of vector bundles on ${\bf P}$:
$$
0\rightarrow
{\cal O}_{\bf P}(-r-1)\otimes_{\bf C}\bigwedge^{r+1}V\rightarrow
{\cal O}_{\bf P}(-r)\otimes_{\bf C}\bigwedge^rV\rightarrow ...
\rightarrow {\cal O}_{\bf P}(-1)\otimes_{\bf C}V\rightarrow
{\cal O}_{\bf P}\rightarrow 0, 
$$
to obtain, for each integer $l\in {\bf Z}$, a double complex of vector
bundles over ${\bf P}$ with general term
$$
{\cal K}^{s,q}={\cal O}_{\bf P}(s)\otimes B_q\otimes
\bigwedge^{l-s-q}V.
$$
Since the rows of this double complex
are exact, it yields a spectral sequence
abutting to zero, with general term
$$
{\cal E}^{s,q}_1={\cal O}_{\bf P}(s)\otimes K_{l-s-q,q}(B,V).
$$
The differential ${\cal O}_{\bf P}\otimes K_{p+1,q}(B,V)
\stackrel{\partial}{\longrightarrow}
{\cal O}_{\bf P}(1)\otimes K_{p,q}(B,V)$, obtained at the first level of
the above spectral sequence by setting
$s=0$ and $p=l-q-1$, is the
morphism we were looking for.
\end{remark}

\noindent
In the real life, we often meet graded
$SV$-modules which satisfy the following property, 
{\bf which will be assumed to be satisfied 
for the rest of the Section},
$$
B_q=0\mbox{ for all } q<0, \mbox{ and } K_{p,0}(B,V)=0\mbox{ for
all } p\geq 1. $$ 
Under these assumptions, we are able
to prove some properties of the
projection of syzygies for $q=1$, as follows.

\bigskip

\begin{lemma}
The map
$
K_{p+1,1}(B,V)\stackrel{H^0(\partial)}{\longrightarrow}
V\otimes K_{p,1}(B,V)
$
is injective.
\end{lemma}
\noindent
{\em Proof.}
We apply $H^0$ to ${\cal K}^{s,q}$.
The new double complex, with general term
$H^0{\cal K}^{s,q}$, gives rise to a spectral sequence such that
$E_1^{s,q}\cong H^0{\cal E}_1^{s,q}$, and $E_\infty^{s,q}=0$
for any pair $(s,q)$, except for $(s,q)=(0,l)$ when
$E_\infty^{0,l}=B_l$. Moreover,  $H^0(\partial)$
coincides to the differential $E_1^{0,1}\rightarrow E_1^{1,1}$.
Because of the assumptions  we have made on $B$,
we get $E_2^{0,1}=E_3^{0,1}=...=E_\infty^{0,1}=0$.
In particular, $\mbox{Ker }H^0(\partial)=0$.

\bigskip

\begin{lemma}
The map $K_{p+1,1}(B,W_x)\stackrel
{\eta_{x}}{\longrightarrow}
K_{p+1,1}(B,V)$ is injective.
\end{lemma}

\noindent
{\em Proof.} Use the short exact sequence
$(*)$ twice.

\bigskip

\begin{remark}
The projection morphism $\pi_x$ is in this case
a genuine corestriction of $\partial _x$.
Moreover, the following relation also
holds: $\partial _x=\big(\mbox{ev}_x\otimes\mbox{id}
_{K_{p,1}(B,V)}\big)
\circ H^0(\partial )$.
\end{remark}

\begin{lemma}
Any nonzero element
$\alpha\in K_{p+1,1}(B,V)$
survives when we project it from a point
$x$ outside a projective subspace of $\bf P$.
\end{lemma}

\noindent
{\em Proof.}
For any nonzero element
$\alpha \in K_{p+1,1}(B,V)$, writing $H^0(\partial )(\alpha )=
\sum v_i\otimes\alpha _i$, with $\alpha _i\in K_{p,1}(B,V)$
linearly independent and $v_i \in V$, we see that
$\partial_x(\alpha )\not= 0$ as long as
the point $x$ does not belong to the
projective subspace of $\bf P$,
$\{ y\in {\bf P},\; v_i\in W_y\mbox{ for all }i\}$.
This subspace is not the whole $\bf P$,
as $H^0(\partial )(\alpha )\not= 0$.


\begin{sectiune}
Projections of syzygies of projective varieties
\end{sectiune}

\noindent
Consider $X\subset {\bf P}$ a nondegenerate
irreducible variety, $x\in {\bf P}$ a point, and $Y\subset {\bf
P}W_x^*$ the image of $X$ by the projection centered in $x$, ${\bf
P}\setminus\{ x\} \longrightarrow {\bf P}W_x^*$. Denote by $I_X$, and
$I_Y$ the homogenous ideals, and by
$S_X$, and $S_Y$ the homogeneous coordinate rings of
$X$, and $Y$ respectively. We have a natural
embedding $S_Y\hookrightarrow S_X$ which induces, for any integer
$p$, an injective map
$
K_{p,1}(S_Y,W_x)\hookrightarrow
K_{p,1}(S_X,W_x).
$
A key fact is that $L_x\otimes K_{p,1}(S_Y,W_x)$ contains all the projected
syzygies, as shown in \cite{Ehb}:

\bigskip

\begin{lemma} {\rm (Ehbauer)}
$\mbox{\rm Im}\big( K_{p+1,1}(S_X,V)
\stackrel{\pi_x}{\longrightarrow}
L_x\otimes K_{p,1}(S_X,W_x)\big)
\subset L_x\otimes K_{p,1}(S_Y,W_x)$, for any $p\geq 1$.
\end{lemma}

\noindent
{\em Proof.} (cf. \cite{Ehb} Section 6)
Since the image of $X$ in
${\bf P}$ is nondegenerate, we have canonical isomorphisms,
for any $p\geq 1$,
$$K_{p+1,1}(S_X,V) \cong K_{p,2}(I_X,V) \cong
\mbox{\rm Ker}\left((I_X)_2\otimes\bigwedge ^pV
\rightarrow (I_X)_3\otimes\bigwedge ^{p-1}V
\right) ;
$$
similarly for $Y$. The statement of Lemma
reduces then to prove that the projection of an
element of $K_{p,2}(I_X,V)$ belongs to
$L_x\otimes K_{p-1,2}(I_Y,W_x)$.

The see this, we choose homogeneous coordinates $X_0,...,X_r$ on ${\bf P}$,
such that $x=[1:0:...:0]$, and all the coordinate points belong to
$X$. An element in $K_{p+1,1}(S_X,V)$ can be seen then as a
collection of quadrics vanishing on $X$, $\big(
Q_{i_1...i_p}\big)_{0\leq i_1<...<i_p\leq r}$, satisfying the
equations: $$\sum_{k\not\in \{k_1,...,k_{p-1}\} } (-1)^{\# \{
k_i<k\} }Q_{k_1...k...k_{p-1}}X_k=0,$$ for all $0\leq
k_1<...<k_{p-1}\leq r$. Projection of syzygies simply means removing
all the $Q_{i_1...i_p}$'s for $i_1\not= 0$,
and renaming $\widetilde{Q}_{j_1...j_{p-1}}
=Q_{0j_1...j_{p-1}}$, for all $1\leq j_1<...<j_{p-1}\leq r$.
A simple analysis shows that $\big(\widetilde{Q}_{j_1...j_{p-1}}
\big)_{1\leq j_1<...<j_{p-1}\leq r}$ actually belongs to
$(I_Y)_2\otimes \bigwedge ^{p-1}W_x$, and the corresponding equations are
still preserved.

\bigskip

In the case of curves, which is the most interesting for our 
purposes, we can actually prove more.

\bigskip

\begin{proposition}
If $L$ is a globally generated line bundle
over the curve $X$, $x\in X$ is a point, 
and $W_x=H^0(X,L-x)$, then, for any $p\geq 1$, 
the natural map $K_{p,1}(X,L-x)\rightarrow K_{p,1}(X,L,W_x)$ 
is injective, and
$\mbox{\rm Im}\big( K_{p+1,1}(X,L)
\stackrel{\pi_x}{\longrightarrow}
L_x\otimes K_{p,1}(X,L,W_x)\big)
\subset L_x\otimes K_{p,1}(X,L-x).$
\end{proposition}
\noindent
{\em Proof.} The injectivity follows from (1.1). 
Let $V=H^0(X,L)$,  $X'$ be the image of $X$ in ${\bf P}V^*$,
and $Y$ be the image of $X'$ in ${\bf P}W_x^*$ under
the projection centered in $x$. 
By means of (1.2), we have $K_{p+1,1}(S_{X'},V) \cong K_{p+1,1}(X,L)$, 
for all $p$. Using (1.1), we can also show that 
$K_{p,1}(S_Y,W_x)\cong K_{p,1}(X,L-x)$. 
We apply (3.1) to conclude.

\bigskip

\begin{corollary}
If $L$ is a globally generated line bundle
over the curve $X$, and
$K_{p+1,1}(X,L)\not= 0$ for an integer $p\geq 1$, then
$K_{p,1}(X,L-x)\not= 0$ for a generic point
$x\in X$.
\end{corollary}

\noindent
{\em Proof.}
Observe that, by means of (2.5), any
nonzero element in $K_{p+1,1}(X,L)\neq 0$
survives in $K_{p,1}(X,L,W_x)$
by projecting from a generic point
$x\in X$, as the image of $X$ in ${\bf P}H^0(X,L)^*$
is nondegenerate.

\begin{sectiune}
Proofs of main results
\end{sectiune}

\noindent
Both Theorem 1, and Theorem 3, follow as immediate consequences 
of a key Lemma which we show next.

\bigskip

\begin{lemma}
Let $X$ be a curve of genus $g$, and $L_0$ be a nonspecial line bundle
on $X$ of degree $d\geq g$. If $L_0$ satisfies $(M_k)$, and
$x_0$ is a point of $X$ such that $L=L_0+x_0$
is globally generated, then $L$ also satisfies $(M_k)$.
\end{lemma}

\noindent
{\em Proof.} 
We want to prove that the vanishing of 
$K_{p,1}(X,L_0)$, for a positive integer $p$,
implies vanishing for $K_{p+1,1}(X,L)$.
We use the following notation:
$$X_1:=\{ x \in X,\;
h^1(X,L-x)=0\} ;$$
it is an open set of $X$, and, since it contains the point $x_0$,
which corresponds to $L_0$, it is moreover nonempty.
Besides, as under vanishing hypotheses the
Koszul complex is actually a complex of vector
bundles over the base of deformation, and
the locus where it is not exact is closed, we remark that, 
for any positive integer $p$, the set
$$
X_p:=\{ x\in X_1,\; K_{p,1}(X,L-x)\not= 0\}
$$
is closed in $X_1$.

Assume $K_{p+1,1}(X,L)\neq 0$. It follows directly from
$(3.3)$ that $K_{p,1}(X,L-x)\not= 0$
for a generic point $x\in X$, hence $X_p$ contains
a nonempty open set of $X_1$. In this case,
it must be the whole $X_1$, and,
since $x_0\in X_1$, we see that $K_{p,1}(X,L_0)\not= 0$.
\\\\
{\bf Proof of Theorem 1.}
Since for any effective divisor
$D$, $L_0+D$ is nonspecial and globally generated,
as soon as $L_0$ itself is nonspecial and globally generated,
an inductive argument allows us to reduce to
the case of a divisor of degree $1$,
which follows from (4.1).
For the last part of the statement,
observe that $H^0(X,L-L_0)\not= 0$,
and thus $L$ itself is of type $L_0+D$.

\bigskip

\begin{corollary}
If $L_0$ is a nonspecial
globally generated line bundle on the curve $X$,
which satifies $(M_k)$, then, for any effective divisor $D$,
any integer $0\leq \gamma <\mbox{\rm deg}(D)$, and
$x_1,...,x_\gamma$ generic points of $X$,
the line bundle $L_0+D-x_1-...-x_\gamma$ also satisfies
the property $(M_k)$.
\end{corollary}

\noindent
{\em Proof.} The set of 
$\gamma$-tuples $(x_1,...,x_\gamma )$, for which
$L_0+D-x_1-...-x_\gamma$ is nonspecial, globally
generated, and satisfies $(M_k)$, is open.
Since $D$ is effective and $0\leq \gamma <\mbox{\rm deg}(D)$,
there exists a set of points of $X$, $\{ x_1,...,x_\gamma \}$,
such that $D-x_1-...-x_\gamma $ is effective, hence
$L_0+D-x_1-...-x_\gamma$
is nonspecial, globally generated and satisfies
the property $(M_k)$. Therefore, the above-mentioned
open set is  nonempty.
\\\\
{\bf Proof of Theorem 3.} 
If $x$ is a point of $X$, and $L_0=K_X+x$,
then $L_0$ is nonspecial, and $H^0(X,K_X)\cong H^0(X,L_0)$.
In particular, for any $p$, $K_{p,1}(X,K_X)\cong K_{p,1}(X,L_0)$,
and thus $L_0$ satisfies $(M_k)$. Then we apply
(4.1) for $L_0$.

\bigskip


\begin{sectiune}
Syzygies of plane curves
\end{sectiune}

\noindent
One of the main results of \cite{Lo} shows
that Green's canonical conjecture is true for a plane curve.
The main idea of his proof was to relate
the Koszul cohomology of the projective plane, 
by means of a long cohomology sequence, to the 
Koszul cohomology of the curve. In this Section, we use a 
similar strategy to prove that the gonality conjecture 
also holds for plane curves.

We begin by pointing out the following
useful fact:

\bigskip

\begin{lemma} For any integers $k\geq 2$, and $p\geq 1$,
$K_{p,1}({\bf P}^2,{\cal O}_{{\bf P}^2}(k))=0$
if and only if $p\geq N(k):=h^0{\cal O}_{{\bf P}^2}(k)-k$.
\end{lemma}

\noindent
{\em Proof.} Let $r=h^0{\cal O}_{{\bf P}^2}(k)-1$.
Green's duality (\cite{Gr1}, Corollary 2.c.10) in this
case translates into
$$K_{p,1}({\bf P}^2,{\cal O}_{{\bf P}^2}(k))\cong
K_{r-p-2,2}({\bf P}^2,{\cal O}_{{\bf P}^2}(-3),
{\cal O}_{{\bf P}^2}(k))^*.$$
For $p\geq N(k)$, we have $r-p-2\leq -3+k$, so
we can apply Theorem 2.2 of \cite{Gr2},
or Theorem 4.1 of \cite{Gr3} to get the
vanishing of $K_{p,1}({\bf P}^2,{\cal O}_{{\bf P}^2}(k))$.
The fact that this bound is sharp easily
follows from \cite{GL1},  by decomposing
${\cal O}_{{\bf P}^2}(k)={\cal O}_{{\bf P}^2}(k-1)
\otimes {\cal O}_{{\bf P}^2}(1)$.

\bigskip

Let $X\subset {\bf P}^2$ be a smooth
plane curve of degree $k+1=d\geq 3$.
We prove:

\bigskip

\begin{proposition}
The pair
$(X,{\cal O}_X(k))$ satisfies the property
$(M_{k-1})$.
\end{proposition}

\noindent
{\em Proof.}
We denote $V=H^0{\cal O}_{{\bf P}^2}(k)=H^0{\cal O}_X(k)$, and
consider the exact sequence of graded
$SV$-modules:
$$
0\longrightarrow \bigoplus
\limits_{q\geq 0}H^0{\cal O}_{{\bf P}^2}(-d+qk)
\longrightarrow \bigoplus
\limits_{q\geq 0}H^0{\cal O}_{{\bf P}^2}(qk)
\longrightarrow \bigoplus\limits_{q\geq 0}H^0{\cal O}_X(qk)
\longrightarrow 0,
$$
together with its associated long cohomology
sequence:
$$
...\rightarrow
K_{p,1}({\bf P}^2,{\cal O}_{{\bf P}^2}(-d),
{\cal O}_{{\bf P}^2}(k))\rightarrow
K_{p,1}({\bf P}^2,{\cal O}_{{\bf P}^2}(k))\rightarrow
K_{p,1}(X,{\cal O}_X(k))\rightarrow
$$
$$
\rightarrow
K_{p-1,2}({\bf P}^2,{\cal O}_{{\bf P}^2}(-d),
{\cal O}_{{\bf P}^2}(k))\rightarrow
...
$$

By means of  Green's Vanishing Theorem (cf. \cite{Gr1},
Theorem 3.a.1), we see that
$K_{p,1}({\bf P}^2,{\cal O}_{{\bf P}^2}(-d),
{\cal O}_{{\bf P}^2}(k))=0$
for all $p\geq 0$, and
$K_{p-1,2}({\bf P}^2,
{\cal O}_{{\bf P}^2}(-d),{\cal O}_{{\bf P}^2}(k))=0$ for
all $p\geq N(k)$.
Therefore, for all $p\geq N(k)$,
we have $K_{p,1}({\bf P}^2,{\cal O}_{{\bf P}^2}(k))
\cong K_{p,1}(X,{\cal O}_X(k))$,
and the Lemma follows as a direct consequence of (5.1).

\bigskip

In this special case, Corollary 2
applied to ${\cal O}_X(k)$ translates
into the following
(the first part of the statement has been known
for a long time: see, for example, \cite{ACGH} p.56)

\bigskip

\begin{corollary}
The curve $X$ is $k$-gonal,  a pencil
of minimal degree being obtained by projecting from
a point of the curve, and the gonality conjecture is
valid for $X$.
\end{corollary}

The property $(M_{k-1})$ is fulfilled for any line bundle ${\cal
O}_X(n)$, with $n\geq k$, and it fails for ${\cal O}_X(k-2)=K_X$.
We can inquire about ${\cal O}_X(k-1)$, and address the question
of whether it satisfies $(M_{k-1})$ or not. The answer is NO:
since the canonical maps $K_{p,1}({\bf P}^2,{\cal O}_{{\bf
P}^2}(k-1))\rightarrow K_{p,1}(X,{\cal O}_X(k-1))$ are injective,
for all positive $p$ (Green's vanishing and the long exact
sequence, as above), the second part of (5.1) applied for $(k-1)$
shows that the property $(M_{k-1})$ fails for $(X,{\cal
O}_X(k-1))$. It is natural then to ask about intermediate line
bundles (see also (6.4) below):
\\\\
{\bf (5.4) Problem.} Find the least $\gamma$ such that,
for any colinear points $x_1,...,x_\gamma \in X$, the property
$(M_{k-1})$ fails for the line bundle ${\cal O}_X(k)\otimes
{\cal O}_X(-x_1-...-x_\gamma )$.


\begin{sectiune}
Syzygies of curves on a Hirzebruch surface
\end{sectiune}

\noindent
Let $\Sigma_e$ be the Hirzebruch surface
of invariant $e$, and denote by $C_0$ the minimal
section of  $\Sigma_e$, and by $f$ a fiber of the ruling.
For any integers $a$, and $b$, we consider the
line bundle on $\Sigma _e$, $H_{a,b}=(a-1)C_0+(b-1)f$,
to whom we attach the integer $N(a,b):=
h^0{\cal O}_{\Sigma_e}(H_{a,b})-a$.
We choose two integers $k\geq 2$, and $m\geq \mbox{max}\{ ke,k+e\}$,
and a curve $X\in |H_{k+1,m+1}|$;
its genus is computed
by the formula $g=(k-1)(m-1-ke/2)$.

\bigskip

We will show that the pair $(X,H_{k,m}|_X)$
satisfies the property $(M_{k-1})$.
For this, we need two preliminary Lemmas.

\bigskip

\begin{lemma}
We have a natural isomorphism
$K_{p,1}(\Sigma_e,H_{k,m})\cong K_{p,1}(X,H_{k,m}|_X)$,
for any $p\geq N(k,m)$.
\end{lemma}

\noindent
{\em Proof.}
We denote by $V=H^0{\cal O}_{\Sigma_e}(H_{k,m})=
H^0(X,{\cal O}_{\Sigma_e}(H_{k,m})|_X)$,
and consider the graded $SV$-modules:
$B'=\bigoplus\limits_{q\geq 0}H^0{\cal O}_{\Sigma_e}(-X+qH_{k,m})$,
$B=\bigoplus\limits_{q\geq 0}H^0{\cal O}_{\Sigma_e}(qH_{k,m})$,
and $A=B/B'$.
We have then a long exact sequence:
$$
...\rightarrow
K_{p,1}(\Sigma_e,-X,H_{k,m})\rightarrow
K_{p,1}(\Sigma_e,H_{k,m})\rightarrow
K_{p,1}(A,V)\rightarrow
$$
$$
\rightarrow
K_{p-1,2}(\Sigma_e,-X,H_{k,m})\rightarrow
...
$$
Since $A_0\cong {\bf C}$,
$A_1\cong H^0(X,{\cal O}_{\Sigma_e}(H_{k,m})|_X)$,
and $A_2\subset H^0(X,{\cal O}_{\Sigma_e}(2H_{k,m})|_X)$,
(1.1) shows that $K_{p,1}(A,V)\cong K_{p,1}(X,{\cal
O}_{\Sigma_e}(H_{k,m})|_X)$.

We apply again Green's Vanishing Theorem
(\cite{Gr1}, Theorem 3.a.1), which implies
$K_{p,1}(\Sigma_e,-X,H_{k,m})=0$ as soon as
$p\geq 0$, and $K_{p-1,2}(\Sigma_e,-X,H_{k,m})=0$
for all $p\geq h^0{\cal O}_{\Sigma_e}(H_{k-1,m-1})+1$.
Since $N(k,m)\geq h^0{\cal O}_{\Sigma_e}(H_{k-1,m-1})+1$.
we get isomorphisms $K_{p,1}(\Sigma_e,H_{k,m})
\cong K_{p,1}(X,H_{k,m}|_X)$, for all $p\geq N(k,m)$, as stated.

\bigskip

\begin{lemma} For all $p\geq N(k,m)$, we have
$K_{p,1}(\Sigma_e,H_{k,m})=0$.
\end{lemma}

\noindent
{\em Proof.} It suffices to prove vanishing for
$p=N(k,m)$ only, so we stick to this case.

We consider $Y$ a curve in the linear system
$|H_{k,m}|$ on $\Sigma_e$. Then $h^0(Y,H_{k,m}|_Y)=
h^0(\Sigma_e,H_{k,m})-1$,
and the surjection $V\rightarrow H^0(Y,H_{k,m}|_Y)$
corresponds to the hyperplane of ${\bf P}V^*$
which cuts out the curve $Y$ on the image of $\Sigma _e$.
Since we work on a rational surface,
(1.3) applies, and  we get an isomorphism
$K_{N(k,m),1}(\Sigma_e,H_{k,m})\cong K_{N(k,m),1}(Y,H_{k,m}|_Y)$.

We prove the Lemma by induction on $k$;
the first step is $k=2$. In this
case, 
$N(2,m)=h^0(Y,H_{2,m}|_Y)-1$;
by the general theory of syzygies
(apply, for example, Theorem 3.c.1 (1) of \cite{Gr1}),
we know that $K_{h^0(Y,H_{2,m}|_Y)-1,1}(Y,H_{2,m}|_Y)=0$.

The induction step: for $k\geq 3$ we know
$K_{N(k-1,m-1),1}(\Sigma_e,H_{k-1,m-1})=0$, and
we wish to prove that  $K_{N(k,m),1}(\Sigma_e,H_{k,m})=0$,
i.e. $K_{N(k,m),1}(Y,H_{k,m}|_Y)=0$.

Lemma (6.1) applied for $Y$
yields  an isomorphism $K_{N(k-1,m-1),1}(Y,H_{k-1,m-1}|_Y)\cong
K_{N(k-1,m-1),1}(\Sigma_e,H_{k-1,m-1})$, which shows that
$K_{N(k-1,m-1),1}(Y,H_{k-1,m-1}|_Y)=0$.
The induction step is completed by means of
the following facts: 
$$H_{k-1,m-1}.Y\geq 2g(Y)+1,$$
$$|(H_{k,m}-H_{k-1,m-1})|_Y|\not= \emptyset,$$ so there exists
$D\in |(H_{k,m}-H_{k-1,m-1})|_Y|$,
$$N(k,m)-N(k-1,m-1)=
(H_{k,m}-H_{k-1,m-1}).H_{k,m}=\mbox{deg}(D),$$
which altogether permit us to apply Theorem
1 for $Y$, and the line bundles $L_0=H_{k-1,m-1}|_Y$,
$L=H_{k,m}|_Y=L_0+D$.

\bigskip

The obvious inequality $H_{k,m}.X\geq 2g+1$, together
with the existence of a $g^1_k$ given by the ruling,
show that the gonality conjecture is verified for $X$
(the fact that $X$ is $k$-gonal was previously proved
in \cite{Ma}, by using completely different methods):

\bigskip

\begin{theorem}
Let $e\geq 0$, $k\geq 2$, and $m\geq\mbox{\rm max}
\{ke,k+e\}$ be three integers, and
$X\equiv kC_0+mf$ be a curve
on the Hirzebruch surface $\Sigma _e$.
Then $X$ is $k$-gonal, a pencil of minimal
degree being given by the ruling, and for
any $L\in\mbox{\rm Pic}(X)$
with 
$h^0(X,L\otimes {\cal O}_{\Sigma _e}(-H_{k,m})|_X)\neq 0$
the property $(M_{k-1})$ holds for $(X,L)$.
\end{theorem}

In the particular case $e=1$, $m=k+1$,
the curve $X$ is the strict transform
of a plane curve $X_0$ of degree $k+1$,
which passes through the center $x_0$ of the
blowup. Via the natural isomorphism
between $X$ and $X_0$, the line bundle
$H_{k,k+1}|_X$ corresponds to
${\cal O}_{X_0}(k)\otimes {\cal O}_{X_0}(-x_0)$.
This shows that (6.3) does a bit better than (5.2).

\begin{sectiune}
The gonality conjecture for generic nodal curves on 
${\bf P}^1\times {\bf P}^1$
\end{sectiune}

\noindent
The aim of this Section is to prove that the gonality conjecture is
valid for a nodal curve on ${\bf P}^1\times {\bf P}^1$
whose singular points are in general position, fact
which makes the proof of Theorem 4 be straightforward;
this idea was inspired by the work \cite{Sch2}.
We freely use the notation of the previous Sections.

We consider three integers $k\geq 3$, $m\geq k$,
$0\leq \gamma \leq k-2$,
$Y\in |H_{k,m}|$ a smooth curve of genus $(k-2)(m-2)$
on ${\bf P}^1\times {\bf P}^1$, $L_0=H_{k-1,m-1}|_Y$,
and $L=H_{k,m}|_Y$. 
Then $(6.3)$ applies for $Y$,
hence the pair $(Y,L_0)$ satisfies $(M_{k-2})$.
In this case, for a set of general points
$\{ x_1,...,x_\gamma\}\subset Y$ the line bundle
$L-x_1-...-x_\gamma$ satisfies the
property $(M_{k-2})$ (cf. (4.2)). Without loss
of generality, this set can be chosen such that
the points $\{ x_1,...,x_\gamma\}$
are in general position in ${\bf P}^1\times {\bf P}^1$,
in the sense that any two of them are not colinear.
We consider furthermore $\Sigma \stackrel{\sigma}
{\rightarrow}{\bf P}^1\times {\bf P}^1$
the blowup of the points $\{ x_1,...,x_\gamma\}$,
$E$ the exceptional divisor, $H=\sigma^*(k-1,m-1)-E$,
and $\widetilde{Y}\in |H|$ the strict transform of $Y$.

Proposition 1 of \cite{Sch2} ensures the
existence of a smooth connected curve in the linear
system $|\sigma^*(k,m)-2E|$; let us denote it
by $X$. The curve $X$ is $k$-gonal (see below), its genus
equals $g=(k-1)(m-1)-\gamma$, and
its projection on ${\bf P}^1\times {\bf P}^1$
is an irreducible curve of type $(k,m)$,  with
assigned ordinary nodes at $\{ x_1,...,x_\gamma\}$.

\bigskip

\begin{proposition}
The pair $(X,H|_X)$ satisfies the property $(M_{k-1})$.
\end{proposition}

\noindent
{\em Proof.}
It is easy to see that 
$h^i{\cal O}_\Sigma (H-X)=0$ for all $i$. Then we denote
$V=H^0{\cal O}_\Sigma (H)\cong H^0(X,{\cal O}_\Sigma (H)|_X)$,
and compute $\mbox{dim}(V)-k=k(m-1)-\gamma$.
In a similar way as in the Proof of (6.1),
we can check that $\mbox{dim}(V)-k\geq h^0
{\cal O}_\Sigma (2H-X)+1$, hence $K_{p,1}(\Sigma ,H)
\cong K_{p,1}(X,H|_X)$ for all $p\geq \mbox{dim}(V)-k$.
Moreover, (1.3) applied to this case shows that
$K_{p,1}(\Sigma ,H)\cong K_{p,1}
(\widetilde{Y},H|_{\widetilde{Y}})$.
Now, the restriction of $\sigma$ to
$\widetilde{Y}$ gives a natural isomorphism
to $Y$; the line bundle $H|_{\widetilde{Y}}$
on $\widetilde{Y}$ corresponds to
the line bundle $L-x_1-...-x_\gamma$
on $Y$. Besides, Riemann-Roch implies
$h^0(Y,L-x_1-...-x_\gamma )-(k-1)=
\mbox{dim}(V)-k$. Since $L-x_1-...-x_\gamma$
satisfies $(M_{k-2})$, i.e. $K_{p,1}
(Y,L-x_1-...-x_\gamma )=0$ for
$p\geq h^0(Y,L-x_1-...-x_\gamma )-(k-1)$, we are done.

\bigskip

Consequently, as $H.X\geq 2g+1$, and $X$ carries
a $g^1_k$, the gonality conjecture is verified for $X$,
and thus Theorem 4 is true.

\newpage

\begin{sectiune}
Application: Curves on Hirzebruch surfaces satisfy Green's
canonical conjecture
\end{sectiune}

\noindent 
In his paper \cite{Gr1}, Green conjectured that, given $X$ a
curve of genus $g$, and denoting by $c$ its Clifford index (we refer to
\cite{Ma2} for a precise definition) all $K_{p,1}(X,K_X)$
vanish for $p\geq g-c-1$. In other words, the line bundle $K_X$ satisfies
the property $(M_c)$.

The aim of this Section is to verify this conjecture for curves on
Hirzebruch surfaces (we use the same notation as in Section 6).

\bigskip

\begin{theorem}
Let $\Sigma_e$ be the Hirzebruch surface of invariant $e$, and $X$ be
a curve on $\Sigma_e$, numerically equvalent to $kC_0+mf$, with
$k\geq 3$ and $m\geq max\{ ke+1, k+1, k+2e
\}$. Then the Clifford dimension of $X$ equals one, 
and Green's canonical conjecture is valid for $X$.
\end{theorem}

\noindent
{\em Proof.}
We start with the exact sequence,
for any $q$:
$$
0\longrightarrow {\cal O}_{\Sigma_e}(qH-X)\longrightarrow
{\cal O}_{\Sigma_e}(qH)\longrightarrow K_X\longrightarrow 0,
$$
where $H=K_{\Sigma_e}+X$.
Observe that $H^0{\cal O}_{\Sigma_e}(H-X)=H^1{\cal O}_{\Sigma_e}
(H-X)=0$, and thus $H^0{\cal O}_{\Sigma_e}(H)\cong H^0(X,K_X)$.
Obviously, their common dimension equals the genus of $X$.
We get next a long exact sequence (see (1.1), and \cite{Gr1} 1.d.4):
$$
...\rightarrow K_{p,1}(\Sigma_e,-X,H)\rightarrow K_{p,1}(\Sigma_e,H)
\rightarrow K_{p,1}(X,K_X)\rightarrow
K_{p-1,2}(\Sigma_e,-X,H)\rightarrow...
$$
By means of Green's vanishing theorem (\cite{Gr1} 3.a.1),
$K_{p-1,2}(\Sigma_e,-X,H)=0$ for 
$p\geq h^0{\cal O}_{\Sigma_e}(2H-X)+1$.
From the inequality $m\geq {\rm max}\{k+2e,ke+1\}$, we see that Lemma
(6.2) applies, and thus
$K_{p,1}(\Sigma_e,H)=0$ for all $p\geq h^0{\cal O}_{\Sigma_e}(H)
-(k-1)=g-k+1$.
To conclude the proof, we remark that, on the one hand,
$h^0{\cal O}_{\Sigma_e}(2H-X)\leq g-k$, 
which implies the vanishing of
$K_{p,1}(X,K_X)$ for all $p\geq g-k+1$, and, on the other hand,
$X$ naturally carries a $g^1_k$.

\bigskip

\begin{remark}
The role of the vanishing result (6.2) is similar
to the one of Green's vanishing 
\cite{Gr2} 2.2 for the case of plane curves, studied
in \cite{Lo}.
\end{remark}

\begin{remark}
By following a similar strategy, one can easily prove
that Green's canonical conjecture is true for certain nodal curves
on ${\bf P}^1\times {\bf P}^1$. To do that, we simply apply (7.1), and
addapt the proof of (8.1) to the new framework. This also 
provides a new proof of 
the main result of \cite{Sch}, which says that Green's canonical conjecture is
true for generic curves whose genera are sufficiently large compared to the
gonality.
\end{remark}

\bigskip

\noindent 
\small {\em Acknowledgements.} I would like to express
my gratitude to F.-O. Schreyer for guiding my first steps in this
topic, as well as for the many enlightening discussions we had,
and suggestions which I benefited from. I am also grateful to C.
Voisin for valuable comments on an early version of
the manuscript.

This work has been carried out while I was
visiting the University of Bayreuth, the
Abdus Salam ICTP (Trieste), and the
Fourier Institute (Grenoble), to
whom I address my thanks for hospitality.

\end{document}